\documentclass[11pt]{amsart}
\usepackage{amsmath}
\usepackage{amssymb,bbm}
\usepackage{epsfig}
\usepackage{a4wide}
\usepackage{pifont}

\newcommand{\N}{\ensuremath{\mathbb{N}^{}_0}}
\newcommand{\Z}{\ensuremath{\mathbb{Z}}}
\newcommand{\R}{\ensuremath{\mathbb{R}}}
\newcommand{\Q}{\ensuremath{\mathbb{Q}}}
\newcommand{\C}{\ensuremath{\mathbb{C}}}
\newcommand{\Sb}{\ensuremath{\mathbb{S}}}
\newcommand{\X}{\ensuremath{\mathbb{X}}}
\newcommand{\Sc}{\ensuremath{{\mathcal S}}}
\newcommand{\D}{\ensuremath{{\mathcal D}}}
\newcommand{\F}{\ensuremath{{\mathcal F}}}
\newcommand{\gL}{\varLambda}
\newcommand{\gG}{\varGamma}
\newcommand{\gT}{\varTheta}
\newcommand{\smc}{\scriptscriptstyle\square}

\DeclareMathOperator{\dd}{d\!}
\DeclareMathOperator{\dens}{dens}
\DeclareMathOperator{\vol}{vol}

\newtheorem{thm}{Theorem}
\newtheorem{prop}[thm]{Proposition}
\newtheorem{lem}[thm]{Lemma}
\newtheorem{cor}[thm]{Corollary}

\newtheorem{claim}{Claim}
\newtheorem{conj}{Observation}

\theoremstyle{definition}

\begin{document}

\title[A radial analogue of Poisson's summation formula]{A radial
  analogue of Poisson's summation formula\\[2mm] with applications to
  powder diffraction\\[2mm] and pinwheel patterns}  

\author{Michael Baake}
\address{Fakult\"at f\"ur Mathematik, Universit\"at Bielefeld, 
Postfach 100131, 33501 Bielefeld,  Germany}
\email{\texttt{mbaake@math.uni-bielefeld.de}, 
\texttt{dirk.frettloeh@math.uni-bielefeld.de}}
\urladdr{\texttt{http://www.math.uni-bielefeld.de/baake},
\texttt{http://www.math.uni-bielefeld.de/baake/frettloe/}}

\author{Dirk Frettl\"oh}

\author{Uwe Grimm}
\address{
Department of Mathematics, The Open University, Walton Hall,
Milton Keynes MK7 6AA, UK}
\email{u.g.grimm@open.ac.uk}
\urladdr{http://mcs.open.ac.uk/ugg2/}

\begin{abstract} 
Diffraction images with continuous rotation symmetry arise from
amorphous systems, but also from regular crystals when investigated
by powder diffraction. On the theoretical side, pinwheel patterns and
their higher dimensional generalisations display such symmetries as
well, in spite of being perfectly ordered. We present first steps and
results towards a general frame to investigate such systems, with
emphasis on statistical properties that are helpful to understand and
compare the diffraction images. An alternative substitution rule for
the pinwheel tiling, with two different prototiles, permits the
derivation of several combinatorial and spectral properties of this
still somewhat enigmatic example. These results are compared with 
properties of the square lattice and its powder diffraction.

\end{abstract} 

\maketitle

\section{Introduction} \label{seq:intro}

Since the discovery of quasicrystals some 20 years ago, mathematicians
and physicists have gained a reasonable understanding of aperiodically
ordered systems, in particular of those obtained from the projection
method. Such sets are called cut and project sets, or model sets
\cite{moo}. These are Delone sets of finite local complexity with
respect to translations, which also means that any finite patch occurs
in finitely many orientations only.

Much less is known about aperiodically ordered systems with local
patches occurring in infinitely many orientations, such as the
pinwheel tiling of the plane or its three-dimensional counterpart,
see \cite{rad} and references therein. Arguably, these are closer to
amorphous systems, but still perfectly ordered. To our knowledge, no
mathematically satisfactory frame for the analysis of radially
symmetric systems and the comparison of their spectral properties has
been developed so far. It is the aim of this contribution to show
first steps in this direction, by combining results and methods from
discrete geometry with the more recent approach of mathematical
diffraction theory.

Our guiding examples are the square lattice and the pinwheel tiling of
the plane. The diffraction of the square lattice is well understood,
and it is not difficult to get some insight into its powder
diffraction. The latter emerges from the presence of many grains in
random mutual orientations, and thus requires a setting with circular
symmetry.

Circular symmetry is also a fundamental property of the pinwheel
tiling, resp.\ its compact hull. This tiling has recently been
reinvestigated from the autocorrelation and diffraction point of view
\cite{mps}. However, hardly any explicit calculation exists in the
literature, the reason being the enigmatic nature of the substitution
generated pinwheel tiling.  New insight is gained by means of an
alternative construction on the basis of a substitution rule with two
distinct prototiles. This gives access to quantities such as
frequencies (hence also to the frequency module), distance sets and
the ring structure of the diffraction measure.

In Section~\ref{sec:rpsf}, we derive a radial analogue of Poisson's
summation formula for tempered distributions, which is needed in the
following analysis. In the last section, we outline the general
structure by presenting a number of results. They are clearly
distinguished according to their present status into theorems (with
proofs or references), claims (with a sketch of the idea and a
reference to future work) and observations (based on numerical or
preliminary evidence). Since our alternative substitution will no
doubt enable other developments as well, we hope that further progress
is stimulated by the results presented in this paper.

\subsection{Notation and preliminaries}

A rotation through $\alpha$ about the origin is denoted by
$R^{}_{\alpha}$, and $B^{}_r(x)$ is the closed ball of radius $r$
centred in $x$.  Whenever we speak of an \emph{absolute frequency}
(e.g., of a point set), we mean the average number per unit volume,
whereas \emph{relative frequency} of a subset of objects is used with
respect to the entire number of objects. The absolute frequency of the
points in a point set $X$ (if it exists) is also called the
\emph{density} of $X$, denoted by $\dens(X)$. The set of non-negative
integers is called $\N$, and the unit circle is
$\Sb^{1}=\{z\in\C:|z|=1\}$.  If $M$ is a locally finite point set, the
corresponding \emph{Dirac comb} $\delta^{}_M := \sum_{x \in M}
\delta^{}_x$, where $\delta^{}_x$ is the normalised point (or Dirac)
measure in $x$, is a well-defined measure. The Fourier transform of
$g$ is denoted by $\widehat{g}$. Whenever we use Fourier transform for
measures below, we are working in the framework of tempered
distributions \cite{sch}.

\section{A radial analogue of Poisson's summation formula}
\label{sec:rpsf}

The diffraction pattern of an ideal crystal, supported on a point
lattice $\gG \subset \R^d$, can be obtained by the Poisson
summation formula 
(PSF) for lattice measures or Dirac combs \cite{sch,cor1,cor2} 
\begin{equation} \label{eq:psf}
\widehat{\delta}^{}_{\gG} = \dens(\gG) \cdot
\delta^{}_{\gG^{\ast}},
\end{equation}
where $\gG^{\ast} := \{ x \in \R^d : x \!\cdot\! y \in \Z$ for all $y
\in \gG \}$ is the dual lattice. The distribution-valued (or
measure-valued) version follows from the classical PSF, compare
\cite{cor1,ik}, via applying it to a (compactly supported) Schwartz
function. It is often used to derive the diffraction measure
$\widehat{\gamma}^{}_{\omega}$ of a lattice periodic measure $\omega =
\varrho \ast \delta^{}_{\gG}$ (with $\varrho$ a finite measure). One
obtains the autocorrelation measure
\begin{equation}  \label{eq:autocor} 
\gamma^{}_{\omega} = \dens(\gG) \cdot (\varrho \ast
\widetilde{\varrho}\, ) \ast \delta^{}_{\gG} ,
\end{equation}
where $\widetilde{\varrho}(g) := \overline{\varrho(\widetilde{g})}$
with complex conjugation $\overline{\vphantom{g}.}$ and
$\widetilde{g}(x)=\overline{g(-x)}$, and the diffraction measure
\[
\widehat{\gamma}^{}_{\omega} = ( \dens(\gG))^2 \cdot |
\widehat{\varrho} \, |^2 \cdot \delta^{}_{\gG^{\ast}} ,
\] 
see \cite{hof,baa} for details. Our interest is to extend this
approach to situations with circular (or spherical) symmetry. 
 
Let us first consider the square lattice $\Z^2$, and recall that its
circular shells have radii precisely in the set 
\begin{equation}\label{defrad}
 \D^{}_{\smc} := \{ r \ge 0 : r^2 = m^2
+ n^2 \; \mbox{with $m,n \in \Z$} \} =
\{0,1,\sqrt{2},2,\sqrt{5},2 \sqrt{2}, 3, \ldots \} . 
\end{equation}
This is the set of non-negative numbers whose squares are integers
that contain primes $p \equiv 3 \mod 4$ only to even
powers. Moreover, on a shell of radius 
$r \in \D^{}_{\smc}$, one 
finds finitely many lattice points, their number being given by 
\begin{equation} \label{etaq}
 \eta^{}_{\smc} (r) =   
\begin{cases} 
1, & r=0, \\
4 a(r^2), & r \in \D^{}_{\smc}\setminus\{0\}. 
\end{cases}
\end{equation}
Here, $a(n)$ is the number of ideals of norm $n$ in $\Z[i]$, the ring
of Gaussian integers. This is a multiplicative arithmetic function,
thus specified completely by its values at prime powers (see
\cite{bg,haw}).  They are given by
\[ 
a(p^{\ell}) = \begin{cases} 1, & \text{if $p=2$}, \\
\ell+1, &  \text{if $p \equiv 1\bmod 4$}, \\ 
0, &  \text{if $p\equiv 3 \bmod 4$ and $\ell$ odd}, \\ 
1, &  \text{if $p\equiv 3 \bmod 4$ and $\ell$ even}. \end{cases}
\] 
As $\Z^2$ is self-dual as a lattice, with $\dens(\Z^2)=1$, the PSF
\eqref{eq:psf} simplifies to $\widehat{\delta}^{}_{\Z^2} =
\delta^{}_{\Z^2}$.

Choose an irrational number $\alpha \in (0,1)$. By Weyl's lemma
\cite{kn}, the sequence $(n \alpha \bmod 1)_{n \ge 1}$ is uniformly
distributed in $(0,1)$. Consider the sequence 
$(\omega^{}_N)^{}_{N \ge  1}$ of measures defined by
\begin{equation} \label{eq:sumrg} 
 \omega^{}_N = \frac{1}{N} \sum_{n=1}^N \delta^{}_{R^n \Z^2}, 
\end{equation}
where $R=R^{}_{2 \pi \alpha}$ is the rotation through $2 \pi \alpha$.
If $|x|=r>0$, the sequence $(R^n x)^{}_{n \ge 1}$ is uniformly
distributed on $\partial B^{}_r(0)$, again by Weyl's lemma. Observe
that all lattices $R^n \Z^2$ share the same set of possible shell
radii, namely $\D^{}_{\smc}$ of \eqref{defrad}.  This implies that,
given an arbitrary compactly supported continuous function $\varphi$,
one has
\[ 
\lim_{N \to \infty} \frac{1}{N} \sum_{n=1}^N \delta^{}_{R^n \Z^2}
(\varphi) = \sum_{r \in \D^{}_{\smc}}
\eta^{}_{\smc} (r) \mu^{}_r (\varphi),
\]  
where the measure $\mu^{}_r$ is the normalised uniform distribution on
$\partial B^{}_r (0) = \{ x \in \R^2 : |x| =r \}$, with
$\mu^{}_0=\delta^{}_0$. This establishes the following result.

\begin{prop}
The sequence $(\omega^{}_N)^{}_{N \ge 1}$ of \eqref{eq:sumrg}
converges in the vague topology, and
\[ 
\omega:= \lim_{N \to \infty} \omega^{}_N  = \sum_{r \in
    \D^{}_{\smc}} 
   \eta^{}_{\smc} (r) \, \mu^{}_r, 
\]
with shelling numbers $\eta^{}_{\smc}(r)$ and
probability measures $\mu^{}_r$ as introduced above.\qed
\end{prop}

It is obvious that the limit is, at the same time, also a limit of
tempered distributions, i.e., a limit in $\Sc'(\R^2)$. As the Fourier
transform is continuous on $\Sc'(\R^2)$, one has
\[ 
\widehat{\omega} = \Bigl( \lim_{N \to \infty} \omega^{}_N
\Bigr)^{\widehat{}}  = \lim_{N \to \infty} \widehat{\omega}^{}_N. 
\] 
Employing the ordinary PSF \eqref{eq:psf}, one finds
\[  
\widehat{\omega}^{}_N = \frac{1}{N} \sum_{n=1}^N
\widehat{\delta}^{}_{R^n \Z^2} = \frac{1}{N} \sum_{n=1}^N
\delta^{}_{(R^n \Z^2)^\ast} = \frac{1}{N} \sum_{n=1}^N \delta^{}_{R^n
\Z^2} = \omega^{}_N 
\] 
where we used that $(R \gG)^{\ast} = R \gG^{\ast}$ for an
isometry $R$, and $\dens(R^n \Z^2)=\dens(\Z^2)=1$.  Combining the last
two equations, one sees that
\[ 
\widehat{\omega} = \lim_{N \to \infty} \widehat{\omega}^{}_N = \lim_{N
  \to \infty} \omega^{}_N = \omega. 
\]

Using polar coordinates, the Fourier transform of $\mu^{}_r$ --- which
is an analytic function because $\mu^{}_r$ has compact support --- can
be expressed by a Bessel function of the first kind via the following
calculation,
\begin{equation}\label{eq:foubessel}
\begin{split} 
\widehat{\mu}^{}_r(k) & = \int_{\R^2} e^{-2\pi i k\cdot x} \dd \mu^{}_r(x) 
 = \int_{0}^{\infty} \frac{1}{2\pi} \int_{0}^{2\pi} e^{-2\pi i
   |k| \varrho \cos{\varphi}} \dd \varphi \dd \delta^{}_r(\varrho)  \\[2mm]
& =  \int_{0}^{\infty} J^{}_0(2\pi|k|\varrho) \dd \delta^{}_r(\varrho) =  
J^{}_0 ( 2\pi |k| r),  
\end{split}
\end{equation}
where $J^{}_0(z) = \sum_{\ell=0}^{\infty} \frac{(-1)^\ell}{(\ell!)^2}
(\frac{z}{2})^{2\ell}$. This yields the identity
\begin{equation} \label{eq:sumbess}
\widehat{\omega} = \sum_{ r \in \D^{}_{\smc}}
\eta^{}_{\smc}(r) J^{}_0(2 \pi |k|r),  
\end{equation}
to be understood in the distribution sense. It has the
following consequence. 

\begin{cor} \label{corsumbess}
With $\eta^{}_{\smc}(r)$ and $\mu^{}_r$ as
introduced above, one has
\[  
    \sum_{r \in \D^{}_{\smc}} \eta^{}_{\smc} (r) \, \mu^{}_r =
\Bigl(  \sum_{r \in \D^{}_{\smc}}
\eta^{}_{\smc} (r) \, \mu^{}_r \Bigr)^{\widehat{}} =
\sum_{r \in \D^{}_{\smc}}
\eta^{}_{\smc} (r) J^{}_0(2 \pi |k| r ),
\] 
where the last expression is to be understood in the distribution
sense. \qed
\end{cor}

Identities of this type can be viewed as measure-valued
generalisations of classic Hardy-Landau-Voronoi summation formulae,
see \cite[Sec.\ 4.4]{ik} and references given there for
details. However, in view of the rather delicate convergence
properties, a direct verification via the PSF \eqref{eq:psf} seems a
simpler approach. 

Observe that one can alternatively reduce the problem to one dimension
and employ a Hankel transform, compare \cite[Sec.\ 4.4]{ik}. This
requires a separate treatment of the transformed measure at the origin
in $k$-space. As this looks slightly artificial from the point of view
of diffraction, we stick to ordinary Fourier transform.

The version for a general lattice $\gG \subset \R^d$ reads as
follows, where, in analogy to $\eta^{}_{\smc}(r)$, 
$\eta^{}_{\gG}(r)$ and $\eta^{}_{\gG^{\ast}}(r)$ denote
the number of points of $\gG$ and $\gG^{\ast}$ on centred
shells $\partial B^{}_r(0)$ of radius $r$. 

\begin{thm}[Radial PSF] \label{radpsf}
Let $\gG$ be a lattice of full rank in $\R^d$, with dual lattice
$\gG^{\ast}$. If the sets of radii for non-empty shells are
$\D^{}_{\gG}$ and $\D^{}_{\gG^{\ast}}$, with shelling
numbers $\eta^{}_{\gG}(r)$ and $\eta^{}_{\gG^{\ast}}(r)$,
the classical PSF \eqref{eq:psf} has the radial
analogue  
\begin{equation} \label{eq:radpsf}
\Bigl( \sum_{ r \in \D^{}_{\gG}}  \eta^{}_{\gG}(r) \, \mu^{}_r
\Bigr)^{\widehat{}} 
= \sum_{r \in \D^{}_{\gG} } \eta^{}_{\gG}(r)
  \, \widehat{\mu}^{}_r
= \dens(\gG) \sum_{r \in  \D^{}_{\gG^{\ast}} }
\eta^{}_{\gG^{\ast}}(r) \, \mu^{}_r ,
\end{equation}     
where $\mu^{}_r$ denotes the uniform probability measure on the sphere of
radius $r$ around the origin. 
\end{thm}

\begin{proof}
Select a sequence of isometries $(R^{}_n)^{}_{n \ge 0}, \, R^{}_n \in \mathrm{SO}(d)$,
such that $(R^{}_n x)^{}_{n \ge 0}$, for a fixed $x$ of length 1, is
uniformly distributed on the unit sphere ${\mathbb S}^d$. Consider
then the sequence of measures defined by 
\[ \omega^{}_N = \frac{1}{N} \sum_{n=1}^N \delta^{}_{R^{}_n \gG} . \]
The claim now follows from Weyl's lemma and the classical PSF
\eqref{eq:psf} in complete analogy to our previous planar example. 
\end{proof}

\begin{figure} 
\epsfig{file=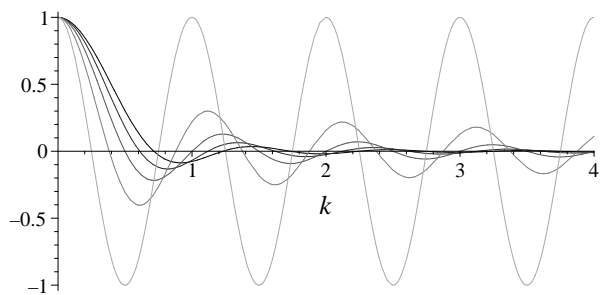, height=40mm}
\caption{\label{j123}
A plot of the radial structure of the function $\widehat{\mu}^{}_1(k)$ 
of \eqref{ftmur} for dimensions $d=1$ (light grey), $2$, $3$, $4$, and $5$ (black). }
\end{figure}

The formula for general $d$ can also be expressed in terms of Bessel
functions of the first kind. Here, by standard calculations with
spherical coordinates and integral representations of Bessel
functions, one obtains
\begin{equation} \label{ftmur}
\widehat{\mu^{}}_r(k) = \int_{\R^d} e^{-2 \pi i kx} \dd \mu^{}_r(x) =
\Gamma\bigl(\frac{d}{2}\bigr)\, \frac{ J^{}_{\frac{d}{2}-1} (2 \pi |k|r)}{(\pi
  |k|r)^{\frac{d}{2} -1}},  
\end{equation}
where $\Gamma(x)$ is the gamma function and $J^{}_{\nu}(z) = (
\frac{z}{2})^{\nu} \sum\limits_{\ell=0}^{\infty} 
\frac{(-1)^\ell}{\ell! \Gamma(\nu + \ell+1)} (\frac{z}{2})^{2\ell}$.
In particular, 
\[ 
\widehat{\mu}^{}_r(k) = \begin{cases} 
\cos(2\pi k r), & \text{if $d=1$}, \\[1mm]
J^{}_0(2 \pi |k| r), & \text{if $d=2$}, \\[1mm]
\frac{\sin(2 \pi |k| r)}{2 \pi |k| r}, & \text{if $d=3$}. 
\end{cases} 
\]
The analogue of Corollary~\ref{corsumbess} for $d=1$ and 
$\gG = \Z$ thus reads 
\[ \sum_{m \in \Z} \cos(2 \pi km) = \widehat{\delta}^{}_{\Z} =
\delta^{}_{\Z}, \]
which is just another form of the ordinary PSF in this case (as radial
averaging is trivial here). Figure~\ref{j123} shows \eqref{ftmur} for
$r=1$ and various values of $d$. 

\section{A simplistic approach to powder diffraction} 

The intensity distribution in powder diffraction emerges from a
collection of grains in random and mutually uncorrelated
orientations. Its precise theoretical description is difficult,
compare \cite{war} and references therein.

Here, we look into a rather simplistic approach that nevertheless
captures the essence of the diffraction image. For simplicity, and for
comparison with related pinwheel patterns, we explain this for the
square lattice $\Z^2$. Instead of working with grains of finite size,
we consider the superposition of entire lattices, with appropriate
weights. Moreover, we make the restriction that there is a common
rotation centre for all lattices, which we choose to be the origin. As
a first step, let us take a look at $\Z^2 \cup R \Z^2$, where $R \in
\mathrm{SO}(2)$ is a generic rotation (by which we mean that it is
\emph{not} an element of the group of coincidence rotations
$\mathrm{SOC}(\Z^2) = \mathrm{SO}(2,\Q)$, see \cite{bg} for more on
this concept).

\begin{lem} \label{diff2gitter}
Let $R \in \mathrm{SO}(2)$ be a generic rotation, so that $\Z^2 \cap R\Z^2 =
\{0\}$. Then, the autocorrelation of\/ 
$\omega = \frac{1}{2} (\delta^{}_{\Z^2} + \delta^{}_{R \Z^2})$ is
\begin{equation} \label {eq:autocor2}
\gamma^{}_{\omega} = \frac{1}{4} \delta^{}_{\Z^2} + \frac{1}{4}
\delta^{}_{R\Z^2} + \frac{1}{2} \lambda, 
\end{equation}
where $\lambda$ is the Lebesgue measure in $\R^2$. 
The diffraction measure of $\omega$ is 
\begin{equation} \label {eq:autocor3}
\widehat{\gamma}^{}_{\omega} = \sum_{x \in \Z^2 \cup R\Z^2} I(x)
\delta^{}_x 
\end{equation}   
with $I(0)=1$ and $I(x)= \frac{1}{4}$ for all \ $0 \ne x \in \Z^2 \cup R
\Z^2$. 
\end{lem}
\begin{proof} 
The first two terms in \eqref{eq:autocor2} are the autocorrelations of
$\frac{1}{2}\Z^2$ and of $\frac{1}{2}R\Z^2$. The third term originates
from the cross-correlation between them, where we used the identity
\[ 
\lambda = \lim_{r \to \infty} \frac{1}{\pi r^2} \sum_{\substack{
    x \in \Z^2 \cap B^{}_r(0) \\ y \in R\Z^2 \cap B^{}_r(0)}}
\delta^{}_{x-y}, 
\] 
with the limit being taken in the vague topology. This identity can be
derived as follows.  It is easy to see that each square in the square
lattice is hit by $x-y$ exactly once, if $x$ is arbitrary but fixed
and $y$ runs through $R\Z^2$. (The exception is $y=0$, but we can
neglect it since it plays no role in the limit.)  Thus each square in
$\Z^2$ is hit with the same frequency in the limit. Since $R$ is a
generic rotation, the sequence $(x-y \bmod (1,1))$ is uniformly
distributed in the fundamental domain $[0,1)^2$ of $\Z^2$, if $x$ is
arbitrary but fixed and $y$ runs through $R \Z^2$. This is a
consequence of the uniform distribution of $(n\alpha \bmod 1)^{}_{n \ge
1}$ in $(0,1)$ if $\alpha$ is irrational. Thus, the points in the sum
above are uniformly distributed in $\R^2$ in the limit, and the series
converges in the vague topology to the Lebesgue measure.

Now, Equation \eqref{eq:autocor3} follows immediately from
$\widehat{\lambda} = \delta^{}_0$ and the PSF \eqref{eq:psf}.  
\end{proof}

For completeness, let us briefly comment on the situation that
$R\in\mathrm{SOC}(\Z^2)$. In this case, $\gT:=\Z^2\cap R\Z^2$ is a
sublattice of $\Z^2$ of finite index, which is $1/\dens(\gT)$.  It is
possible to show that $\omega=\frac{1}{2}(\delta^{}_{\Z^2}
+\delta^{}_{R\Z^2})$ has diffraction
\[
\widehat{\gamma}^{}_{\omega} = 
\delta^{}_{\gT} + \frac{1}{4}\delta^{}_{\Z^2\setminus\gT} +  
\frac{1}{4}\delta^{}_{R\Z^2\setminus\gT}.
\]
The difference to the diffraction formula in Lemma~\ref{diff2gitter}
is concentrated on $\gT\setminus\{0\}$, and hence plays no role in our
further discussion when $\dens(\gT)\to 0$, which happens under
multiple coincidence intersections \cite{bg}. We may thus disregard
coincidence rotations for our present purposes.

\begin{figure} 
\epsfig{file=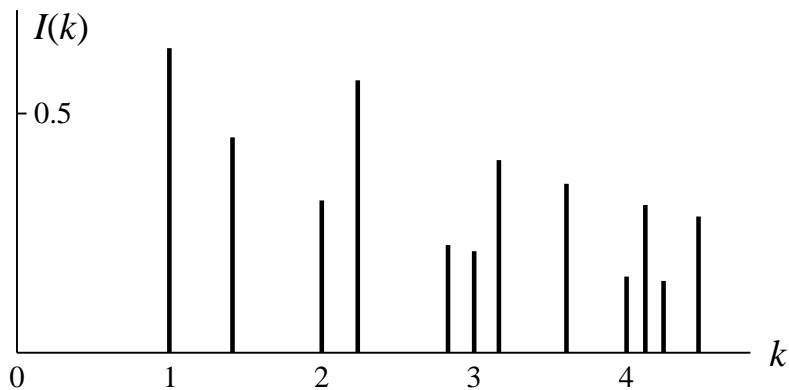,width=0.7\textwidth}
\caption{\label{sqexact}
Radial dependence of ring intensities of square lattice powder diffraction.
The central intensity is not shown, see text for details.}
\end{figure}

\begin{figure} 
\epsfig{file=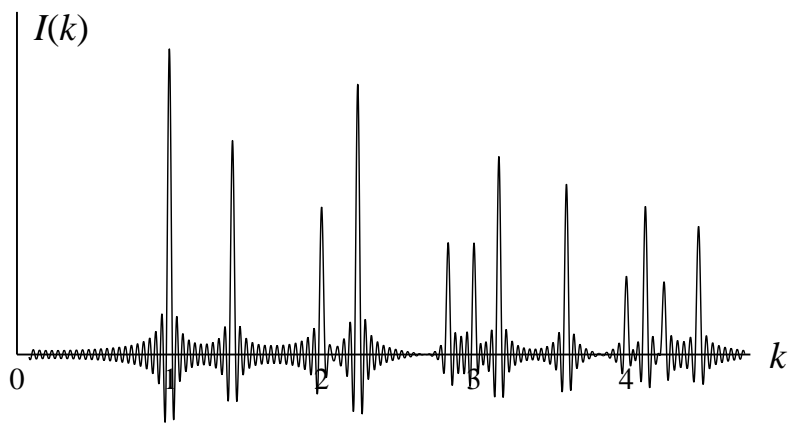,width=0.7\textwidth}
\caption{\label{sqnum}
Numerical approximation of Figure~\ref{sqexact}, 
based on summing Bessel functions, for
radii $r\le 25$. The required radial autocorrelation coefficients
$\eta^{}_{\smc}(r)$ are taken from Equation~\eqref{etaq}.
}
\end{figure}

Let us continue by considering multiple intersections.  If all
rotations are generic (in the sense that the lattices $R_{i}\Z^2$ and
$R_{j}\Z^2$ are distinct apart from the origin) and satisfy a uniform
distribution property, we obtain the following result.

\begin{prop}
Let $\gG = \Z^2$, $\gG^{}_j = R^{}_j \Z^2$ with
generic $R^{}_j$, i.e., $\gG^{}_j \cap \gG^{}_k = \{ 0 \}$
for $j \ne k$, and define $\omega = \frac{1}{N} \sum_{j=1}^N
\delta^{}_{\gG^{}_j}$. 
Then, one has the identities
\[ 
\gamma^{}_{\omega} = 
\frac{N-1}{N} \lambda \; +
\; \sum\limits_{j=1}^N \frac{1}{N^2}
\delta^{}_{\gG^{}_j}  \quad \mbox{ and }
\quad \widehat{\gamma}^{}_{\omega} = \frac{N-1}{N} \delta^{}_0 +
\frac{1}{N} \Bigl( \frac{1}{N} \sum\limits_{j=1}^N
\widehat{\delta}^{}_{\gG^{}_j}   \Bigr). \]
If $( R^{}_j x )^{}_{j \ge 0}$ is uniformly distributed on $\Sb^1$ for
some fixed $x \in \Sb^1$, we obtain 
\[ \lim_{N\to \infty} \frac{1}{N} \sum_{j=1}^N
\widehat{\delta}^{}_{\gG^{}_j} = 
\sum_{r \ge 0} \eta^{}_{\smc}(r) \mu^{}_r,   \] 
with $\eta^{}_{\smc}(r)$ as in \eqref{etaq}.
\end{prop}
 
\begin{proof}
The first statement --- about $\gamma^{}_{\omega}$ --- follows 
as in Lemma~\ref{diff2gitter}.  
Observing $(\gG^{}_j)^{\ast} = R_j \gG^{\ast}$, this implies 
\[ 
\widehat{\gamma}^{}_{\omega} =  \sum_{x \in \bigcup_{j}
  R_j \gG^{\ast}} I(x) \delta^{}_x, 
\]
where $I(0)=1$ and $I(x)=\frac{1}{N^2}$ otherwise. Consequently, one
finds 
\[ 
\widehat{\gamma}^{}_{\omega} = \frac{N-1}{N} \delta^{}_0 +
\frac{1}{N} \Bigl( \frac{1}{N} \sum_{j=1}^N
\widehat{\delta}^{}_{\gG^{}_j} \Bigr). 
\] 
By Theorem~\ref{radpsf} in connection with Weyl's lemma and the
self-duality of $\Z^2$, the term in the brackets converges to $\sum_{r
\ge 0} \eta^{}_{\smc}(r) \mu^{}_r$ as $N \to \infty$.
\end{proof}

This simple argument shows that, after discarding the central
intensity and multiplying the remainder by $N$, one is left with a
circular diffraction pattern in the spirit of the radial PSF in
Theorem~\ref{radpsf}. Disregarding the central intensity, the shelling
numbers $\eta^{}_{\smc}(r)$ reflect the total intensity, integrated
over the rings of radius $r$, in this idealized approach to powder
diffraction. Consequently, in a measurement that displays the
intensity along a given direction, the resulting radial dependence is
given by $\frac{\eta^{}_{\smc}(r)}{2 \pi r}$, compare
Figure~\ref{sqexact}. The numerical approximation in
Figure~\ref{sqnum}, included for later comparison, shows a strong
oscillatory behaviour, as a result of summing Bessel functions.  Note
that this approximation disregards positivity in favour of including
circular symmetry. The comparison gives a good impression on the
overshooting that originates from this approach.

\section{Application to tilings with statistical circular symmetry} 

\subsection{The pinwheel tiling}
\begin{figure}
\epsfig{width=\textwidth,file=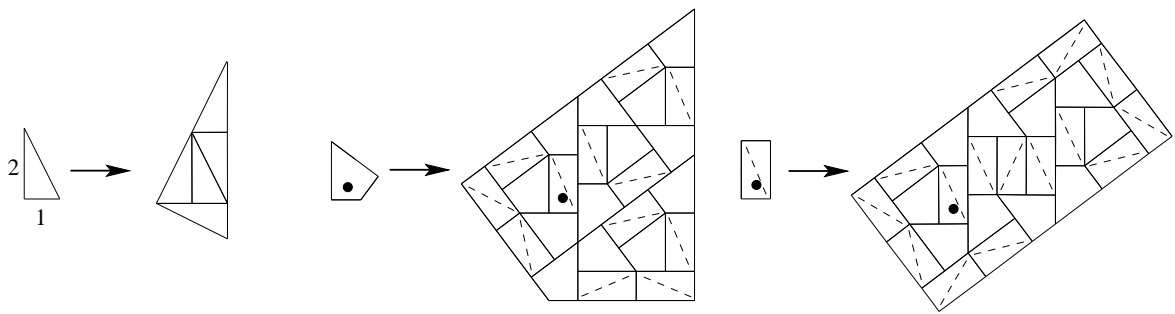}
\caption{The substitution for the pinwheel tiling (left) and for the
kite domino tiling (right). The dashed lines indicate how the dominos
have to be dissected to obtain a pinwheel tiling. They are also needed
to turn the global substitution rule for $KD$ into a local one for 
its hull.
\label{fig:subst}}
\end{figure}

The prototile of the pinwheel tiling is a rectangular triangle with
side length $1$, $2$ and $\sqrt{5}$.  The smallest angle in $T$ is
$\arctan(\frac{1}{2})$.  Here, we choose the prototile $T$ with
vertices $(-\frac{1}{2},-\frac{1}{2})$, $(\frac{1}{2},-\frac{1}{2})$,
$(-\frac{1}{2},\frac{3}{2})$, and equip $T$ with a control point at
$(0,0)$.  Every tile in a pinwheel tiling is either of the form
$R^{}_{\alpha} T + x$, or of the form $R^{}_{\alpha} S T +x$ for some
$x \in \R^2$, where $S$ denotes the reflection in the horizontal axis
and $R^{}_{\alpha}$ rotation through $\alpha$. The substitution
$\sigma$ for the pinwheel tiling is shown in Figure~\ref{fig:subst}
(left). By the action of $\sigma$, $T$ is expanded, rotated by
$-\arctan(\frac{1}{2})$ and dissected into five triangles that are
congruent to $T$. Formally, we define $\sigma(T) := \{ Q_1(T)+x_1,
Q_2(T)+x_2, Q_3(T)+x_3, Q_4(T)+x_4, Q_5(T)+x_5\}$, with $Q_i \in
\mathrm{O}(2)$ and $x_i \in \R^2$. The appropriate choices of $Q_i$
and $x_i$ can be derived from the figure. In particular, $\sigma(T)$
contains a triangle \emph{equal} to $T$, thus we may choose $Q_1 =
\mathbbm{1}$ and $x_1=0$. The substitution $\sigma$ extends in a
natural way to all isometric copies $QT+x$ (where $Q \in
\mathrm{O}(2)$ and $x \in \R^2$) by
$\sigma(QT+x):=Q\sigma(T)+\sqrt{5}x$. Since one of the triangles in
$\sigma(T)$ equals $T$, the tiling $PW$ can be defined as a fixed
point of the substitution $\sigma$. This results in $\sigma(PW)=PW=
\bigcup_{n \ge 0} \sigma^n(T)$. This convention follows \cite{mps} but
deviates from \cite{rad}, because it is advantageous for us to include
the rotation in the definition of $\sigma$.

\subsection{The kite domino tiling}
It is advantageous to consider a second substitution which yields a
closely related tiling. Consider the two prototiles $K$ and $D$ (for
'kite' and 'domino'), where $D$ is the rectangle with vertices
$(-\frac{1}{2},-\frac{1}{2})$, $(\frac{1}{2},-\frac{1}{2})$,
$(\frac{1}{2},\frac{3}{2})$, $(-\frac{1}{2},\frac{3}{2})$, and $K$ is
the quadrangle with vertices $(-\frac{1}{2},-\frac{1}{2})$,
$(\frac{1}{2},-\frac{1}{2})$, $(\frac{11}{10},\frac{3}{10})$,
$(-\frac{1}{2},\frac{3}{2})$. Both tiles consist of two copies of the
pinwheel triangle $T$, glued together along their long edges. Every
pinwheel tiling gives rise to a kite domino tiling by deleting all
edges of length $\sqrt{5}$. The substitution for the kite domino
tiling is shown in Figure~\ref{fig:subst} (right). Essentially, it is
$\sigma^2$ applied to the new prototiles, where $\sigma$ is the
substitution for the pinwheel tiling. If $\varrho$ denotes this new
substitution, a fixed point is given by $KD = \bigcup_{n \ge 0}
\varrho^n(D)$. This is a \emph{global} substitution that works on the
specific $D$ defined above, compare \cite{f} for a discussion of the
general substitution concept. The two tilings, $PW$ and $KD$, are
mutually locally derivable (MLD) in the sense of
\cite{bsj}. Essentially, this means that one is obtained from the
other by local replacement rules and vice versa.  These rules are
evident from Figure~\ref{fig:subst}. 

As $PW$ and $KD$ are MLD, the global substitution $\varrho$ also
induces a \emph{local} substitution in the sense of \cite{bs}, where
two (oriented) copies of the domino have to be distinguished. Loosely
speaking, this follows from the observation that the local surrounding
of any domino in $KD$ or one of the other elements in the hull defined
by $KD$ determines the type of the domino, and hence how to apply
the substitution to it.

We equip $K$ with control points at $(0,0),(\frac{2}{5},\frac{1}{5})$,
and $D$ with control points at $(0,0),(0,1)$. Then, the set of all
control points in $KD$ is equal to the set of control points in
$PW$. This specific discrete point set is a Delone set, denoted by
$\gL$ in the sequel. It is not hard to see that $\gL$ is MLD with both
$PW$ and $KD$, even though one direction of the replacement rule is
less obvious than the one linking $PW$ and $KD$.

Observe that the relative and absolute frequencies of triangles $T$ in
the pinwheel tiling are both equal to $1$. For the relative
frequencies, this is clear since there is only one prototile. The
absolute frequency then follows from the fact that this prototile has
area $1$. Every triangle of $PW$ carries exactly one control point, so
$\dens(\gL)=1,$ too. For the kite domino tiling, standard
Perron-Frobenius theory yields the relative frequencies of kite
(resp.\ domino) as $\frac{5}{11}$ (resp.\ $\frac{6}{11}$). Thus, the
absolute frequency of a kite (resp.\ domino) is $\frac{5}{22}$
(resp.\ $\frac{6}{22}$). The substitution matrix can be extracted from
Figure~\ref{fig:subst}.

\begin{figure}
\epsfig{file=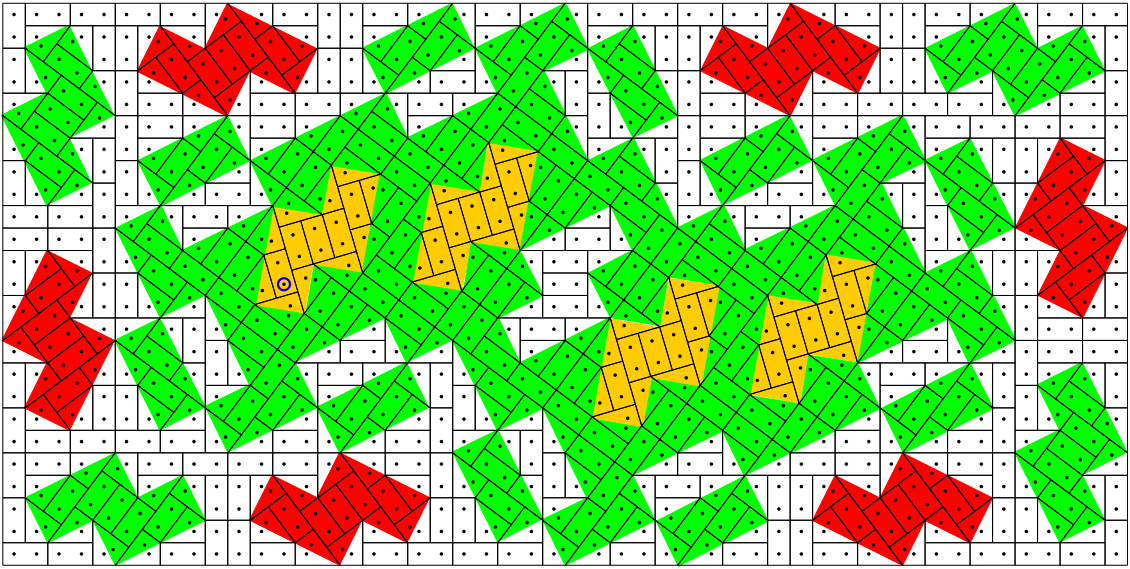, width=150mm}
\caption{A patch from a kite domino tiling. The control points of the
  white tiles are contained in $\Z^2$. The control points in the other
  tiles are contained in rotated copies of $\Z^2$, where the common
  rotation centre is indicated in the figure. The possible rotations
  are described in Claim~\ref{rotgitter}. \label{fig:kdpatch}}
\end{figure}

Discrete structures that are MLD lead to dynamical systems that are
topologically conjugate \cite{kel}. Therefore, we formulate the hull on
the basis of the Delone sets of the control points. Let $\gL$ be one
such set, e.g., the set $\gL$ of control points of $PW$ as defined
above. Define the orbit closure
\[ 
\X(\gL) = \overline{ \R^2 + \gL}^{\mathsf{LRT}}  
\]
in the local rubber topology (LRT) \cite{bl}, which is compact. This
topology is slightly different from the one introduced in
\cite{rad}. But in the present case, both topologies yield the same
hull and are equivalent on $\X(\gL)$.

The following result is well-known \cite{rad,mps}. 

\begin{thm} \label{autoconst}
The hull $\, \X(\gL)$ is \ $\Sb^1$-symmetric. Moreover, $(\X(\gL),
\R^2)$ is a strictly ergodic dynamical system, i.e., it is uniquely
ergodic and minimal.  All elements of\/ $\X(\gL)$ possess the same
autocorrelation measure $\gamma = \gamma^{}_{\gL}$ and the same
diffraction measure $\widehat{\gamma}^{}_{\gL}$. Both measures
are\/ $\Sb^1$-symmetric. \qed
\end{thm}

In this sense, speaking of the diffraction of the pinwheel tiling or
its control points has a unique meaning. At this point, we know that
\[ 
\widehat{\gamma} = \left( \dens(\gL) \right)^2 \delta^{}_0 +
(\widehat{\gamma})^{}_{\mathsf{cont}} = \delta^{}_0 +
(\widehat{\gamma})^{}_{\mathsf{cont}} 
\]
because a translation bounded circularly symmetric measure cannot
contain Bragg peaks other than the trivial one at $k=0$.  That the
intensity coefficient of $\delta^{}_0$ is $\left( \dens(\gL)
\right)^2$, hence $1$ in our case, is a consequence of the equation
\[ 
\widehat{\gamma}(\{0\}) = \lim_{r \to \infty} \frac{1}{\pi r^2}
\, \gamma( B^{}_r (0)) = \left( \lim_{r \to \infty} \frac{1}{\pi r^2}
\, \delta^{}_{\gL} ( B^{}_r (0)) \right)^2 = \left( \dens( \gL )
\right)^2. 
\] 
It can be proved by means of a Fourier concentration argument in
connection with results from \cite{al}; alternatively, see \cite{hof}
for another approach and \cite{mps} for details on the case at hand.
It remains to determine the structure of
$(\widehat{\gamma})^{}_{\mathsf{cont}}$ more closely, in particular
the separation into singular and absolutely continuous components.

\subsection{Results and observations}

The following claims and conjectures are stated for the particular
tiling $KD$ or for its set of control points $\gL$. Since mutual local
derivability extends to entire hulls, we remain in the situation of
Theorem~\ref{autoconst}. Proofs of the claims are either outlined here
or will appear in \cite{bfg}. Let us start with our main observation.

\begin{conj} \label{mainconj}
The diffraction measure $\widehat{\gamma}$ of the pinwheel tiling is of the
form  
\[ 
\widehat{\gamma} = \delta^{}_0 + (\widehat{\gamma})^{}_{\mathsf{sc}} +
(\widehat{\gamma})^{}_{\mathsf{ac}}.
\]
There is a countable set $\D^*\subset\R_{\ge 0}$ such that
\[ 
(\widehat{\gamma})^{}_{\mathsf{sc}} = \sum_{r \in \D^*\setminus\{0\}} 
I(r) \, \mu^{}_r. 
\]
The set $\D^*$ seems to be locally finite, i.e., discrete and closed.
Moreover, $(\widehat{\gamma})^{}_{\mathsf{ac}}$ seems to be 
non-vanishing.
\end{conj}

The set $\D^*$ is a subset of another set, $\D= \D^{}_{\gL}$, which
shows up in the determination of the autocorrelation and is described
in more detail below.

\begin{claim} \label{rotgitter}
Let $\gL$ be the set of all control points of $KD$. Then, $\gL$ is a
Delone subset of a countable union of rotated square lattices. More
precisely, with $\theta = 2 \arctan(\frac{1}{2})$,   
\[ \gL \subset  \bigcup_{k \in \Z}  R^{}_{k \theta}(\Z^2). \]
\end{claim}

Using the geometry of the kite domino tiling, it is not difficult to
establish this property. Figure~\ref{fig:kdpatch} may serve as a
visualisation on a small scale.  This structure suggests that the
diffraction of the pinwheel tiling might share some features with that
of the powder diffraction of $\Z^{2}$, see the comment after the
proof of Lemma~\ref{diff2gitter}. Since $R^{}_{\theta}$ acts as
multiplication by $\frac{1}{5} \bigl( \begin{smallmatrix} 3 & -4 \\ 4
& 3 \end{smallmatrix} \bigr)$, the next result is immediate.  Note
that $R^{}_{\theta}$ is well known from the coincidence site lattice
problem of the square lattice, compare \cite{bg}, and has played a
crucial rule in the explicit calculations in \cite{mps}.
 
\begin{claim} \label{distsupset}
Let $\gL$ be the set of all control points of $KD$. It satisfies
\[ 
\gL \subset \bigl\{ (\textstyle{\frac{n}{5^{k}}, \frac{m}{5^{k}}}) :
m,n \in \Z, \; k\in\N \bigr\}, 
\]  
and the distance set $\D = \D^{}_{\gL} :=\{ |x-y| : x, y \in
\gL \}$ is a subset of  $\{ \sqrt{ \frac{p^2+q^2}{5^{\ell}} } :
p,q,\ell \in \N \}$. It is the same set for all elements of\/
$\X(\gL)$. \qed
\end{claim}

We believe that this subset relation is quite sharp, i.e., that the
difference between the two sets is not too large, in the following sense. 

\begin{claim} \label{allvalues}
All values $r = \sqrt{p^2+q^2}$ with $p,q \in \N$ occur in
$\D$. Moreover, for each $\ell \in \N$, there are infinitely many $p,q
\in \N$ such that $\sqrt{\frac{p^2+q^2}{5^{\ell}}}$ is contained in
$\D$.
\end{claim}

In fact, numerical computations suggest the following, stronger
property. 

\begin{conj}\label{conj2}
For each $\ell \in \N$, $\D$ contains all but finitely many
values of the form $\sqrt{\frac{p^2+q^2}{5^{\ell}}}$.
\end{conj}

In the sequel, the absolute frequencies of distances are of interest,
wherefore we define the \emph{radial autocorrelation function}
\begin{equation} \label{eq:etar}
\eta(r) = \lim_{R \to \infty} \frac{1}{\vol(B^{}_R)} \sum_{
  \substack{ x,y \in \gL\cap B^{}_R (0) \\ |x-y| =r  } } 1. 
\end{equation}
The limit exists due to unique ergodicity, see Theorem
\ref{autoconst}. This permits us to write the autocorrelation of $\gL$
--- and hence of $\X(\gL)$ --- as
\begin{equation} \label{eq:sumrind} 
 \gamma = \sum_{r \in \D} \eta(r) \, \mu^{}_r,   
\end{equation} 
with $\mu^{}_r$ as in Section~\ref{sec:rpsf}. This follows from
Theorem~\ref{autoconst} together with Claim~\ref{allvalues}.  Note
that $\gL$ is repetitive \cite{mps} (defined up to congruence), so
that all $\eta(r)$ in \eqref{eq:sumrind} are strictly positive, by the
definition of $\D$ in Claim~\ref{distsupset}.  Clearly, $\gamma$ is
both a positive and positive definite measure, and a tempered
distribution on $\R^2$.

The following theorem holds in more generality than the context
of this paper. In fact, it holds for all substitution tilings which are
of \emph{finite local complexity} (FLC) and \emph{self-similar}. The
latter means that $\lambda T^{}_i = \bigcup_{T \in \sigma(T^{}_i)} T$. Roughly
speaking, this means that the support of the substitution of each
prototile $T^{}_i$ is similar to $T^{}_i$. Finite local complexity means
that, for some $R>0$, the tiling contains only finitely many local
patches of diameter less than $R$, up to congruence. Note that FLC is
frequently defined with respect to translations, whereas we define FLC
here with respect to congruence, because of the nature of the pinwheel
tiling. Both properties, FLC and self-similarity,
hold for the majority of substitution tilings in the literature
\cite{fh}.

\begin{thm} \label{qlambdad}
In any self-similar substitution tiling of finite local complexity
{\rm (}w.r.t.\ congruence{\rm )} with substitution factor $\lambda$,
all relative frequencies are contained in $\Q(\lambda^d)$.  Moreover,
the tiling can be scaled such that all absolute frequencies are
contained in $\Q(\lambda^d)$ as well.  In particular, if $\lambda^d$
is an integer, and if the tiling is appropriately scaled, all
relative and absolute frequencies are rational.
\end{thm}

For a proof, we refer to \cite{bfg}. 
This theorem, applied to the pinwheel tiling or the kite domino
tiling, yields that all frequencies are rational. 

\begin{claim} All values of $r^2\le 5$ together with 
  $\eta(r)$ are given in the following table. Values marked with an
  asterisk are numerically based conjectures, all
  other values are exact. 
\[ 
\begin{array}{l||c|c|c|c|c|c|c|c|c|c|c|c|c}
\rule[-2mm]{0mm}{7mm} r^2 & 0 & \frac{1}{5} & 1 & \frac{8}{5} &
\frac{9}{5} & \frac{49}{25} & 2 & \frac{13}{5} & \frac{81}{25} &
\frac{17}{5} & 4 & \frac{113}{25} & 5 \\ \hline 
\rule[-2mm]{0mm}{7mm} \eta(r) & 1 & \frac{5}{11} &
\frac{439}{165} & \frac{1}{2} & \frac{67}{165} & \frac{4}{165}
&  \frac{7}{2}^{\ast} & \frac{142}{165} & \frac{4}{165} &
\frac{10}{11}^{\ast} &  3^{\ast} &  \frac{8}{165}^{\ast} &
\frac{73}{15}^{\ast} 
\end{array}
\] 
\end{claim}\smallskip

The value $\eta(0)$ is the absolute frequency of the control points.
Since each triangle in the pinwheel tiling has area $1$, we have
$\eta(0)=1$.  The distance $|x-y|=1/\sqrt{5}$ occurs precisely once
within each kite. Kites have absolute frequency $\frac{5}{22}$, but
the distance $|x-y|$ has to be counted twice in view of formula
\eqref{eq:etar}, hence $\eta(1/\sqrt{5})=\frac{5}{11}$. The other
exact values contained in the table can be established similarly, but
require more sophistication.  With some further effort, one can
determine the \emph{frequency module} of $\gL$, which is the $\Z$-span
of the absolute frequencies of all finite subsets of $\gL$, the latter
standardised to having density $1$ (which is our natural setting here).

\begin{figure}
\epsfig{file=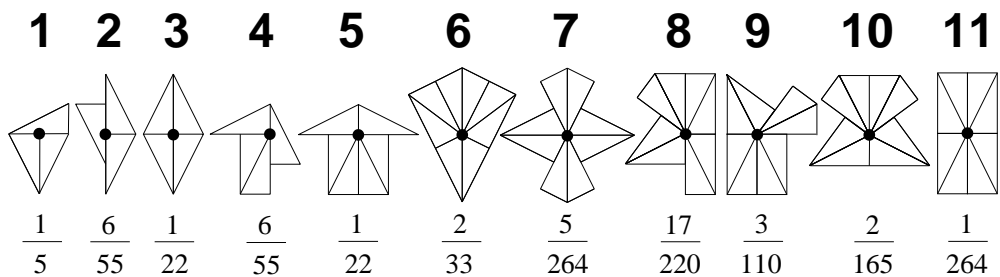}
\caption{The 11 vertex stars of the pinwheel tiling (up to
  congruence) and their absolute frequencies (i.e., frequencies per unit
  area). \label{fig:vertstars}} 
\end{figure}

\begin{claim}
The frequency module of $\gL$, and hence of\/ $\X(\gL)$, is 
$\F= \{ \frac{m}{264 \, \cdot \,  5^{\ell}} \, | \, m \in \Z, \; \ell
\in \N \}$. In particular, it is countably, though not finitely generated. 
\end{claim}

This result is derived from the absolute frequencies of the vertex
stars in the pinwheel tiling, see Figure~\ref{fig:vertstars}. These,
in turn, can be derived from the frequencies of kites and darts in the
kite domino tiling. Details will be given in \cite{bfg}.  

We proceed by considering the Fourier transform in relation to the
distance set $\D$. Since the Fourier transform is continuous
on the space $\Sc'(\R^2)$  of tempered distributions,
\eqref{eq:sumrind} becomes  
\begin{equation}\label{eq:pinbessel}
\widehat{\gamma}^{}_{\gL} = \Bigl( \sum_{ r \in \D}  \eta(r)
\, \mu^{}_r \Bigr)^{\widehat{}} = \sum_{r \in \D }
\eta(r) \, \widehat{\mu}^{}_r, 
\end{equation} 
where the sum is to be understood in the distribution sense. By
Bochner's theorem, $\widehat{\gamma}^{}_{\gL}$ is again a positive
and positive definite measure, and the equation can be understood as a
vague limit as well. 

Let us apply Theorem~\ref{radpsf} to the pinwheel pattern $\gL$.
According to Claim~\ref{distsupset}, the distances in $\D$ arise from
lattices of the form $\frac{1}{5^{\ell/2}}\Z^2$ with $\ell\in\N$.  The
corresponding dual lattices, which are relevant for the diffraction
analysis, are $5^{\ell/2}\Z^2$, which have distance sets
$5^{\ell/2}\D^{}_{\smc}$. Since $\Z^2$ contains a lattice congruent to
$\sqrt{5} \Z^2$, it follows that $5^{\ell / 2} \D^{}_{\smc} \subset
\D^{}_{\smc}$ for all $\ell \in \N$, and that $\D^{}_{\smc}$ contains
all distances that seem relevant from this point of view. At this
stage, we have not found a compelling argument to exclude the
possibility that certain subsets of $\gL$ contribute relevant distances
in a coherent fashion, and thus --  by unique ergodicity -- further
rings $I(r) \mu_r$ (with $r\not\in\D^{}_{\smc}$) to
$\widehat{\gamma}$. However, on the basis of
Theorem~\ref{radpsf}, it is at least plausible that the set $\D^*$ is
closely related with $\D^{}_{\smc}$, or even equal to it. In this
case, $\widehat{\gamma}^{}_{\gL}$ would show rings $I(r) \, \mu^{}_r$
for $r \in \D^{}_{\smc}$ only. According to Claim \ref{allvalues}, all
distances have positive frequencies, so that this assumption would
lead to
\[ 
( \widehat{\gamma})^{}_{\mathsf{sing}} = \sum_{r \in
  \D^{}_{\smc}} I(r) \, \mu^{}_r = \delta^{}_0 +
( \widehat{\gamma})^{}_{\mathsf{sc}} = \delta^{}_0 +  \sum_{0 < r \in
  \D^{}_{\smc}} I(r) \, \mu^{}_r 
\]
with $I(0)=1$ and $I(r)>0$ for all $r \in \D^{}_{\smc}$. 

Let us continue our discussion on the basis of this hypothesis.
If Observation~\ref{conj2} holds, then, for each $\ell \ge 0$, only
finitely many distances of type $\sqrt{\frac{p^2+q^2}{5^{\ell}}}$ are
missing, each one contributing an absolutely continuous part to
$\widehat{\gamma}$. Since $\widehat{\gamma}$ exists as a translation
bounded measure, and since $\widehat{\gamma}^{}_{\mathsf{sing}}$ is
already covered by the results above, every additional contribution
has to contribute to $\widehat{\gamma}^{}_{\mathsf{ac}}$. This
motivates the conjecture that $(\widehat{\gamma})^{}_{\mathsf{ac}} \ne
0$.

\begin{figure} 
\epsfig{file=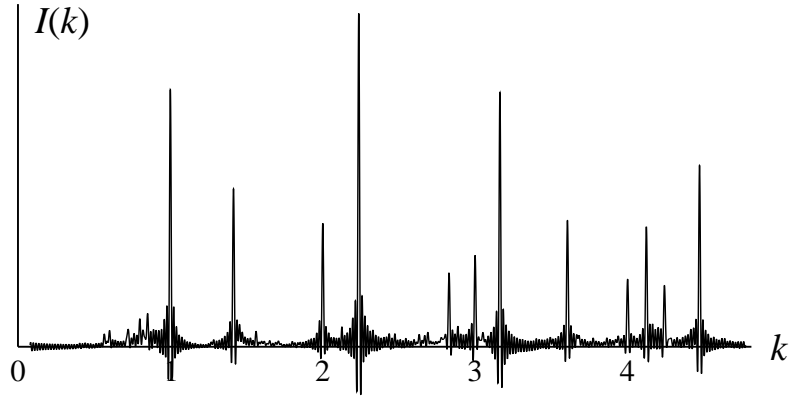,width=0.7\textwidth}
\caption{\label{pin5} Numerical approximation of the radial intensity
dependence of the pinwheel diffraction pattern. It is based on Equations
\eqref{eq:foubessel} and \eqref{eq:pinbessel}, disregarding contributions
to the central intensity. The radial autocorrelation coefficients $\eta(r)$ 
are estimated from a patch of radius of about $56$ (thus effectively using the
fifth iteration of the substitution $\sigma$). The vertical scale is arbitrary
in the sense that it has no meaning without local averaging and integration.}
\end{figure}

\begin{figure} 
\epsfig{file=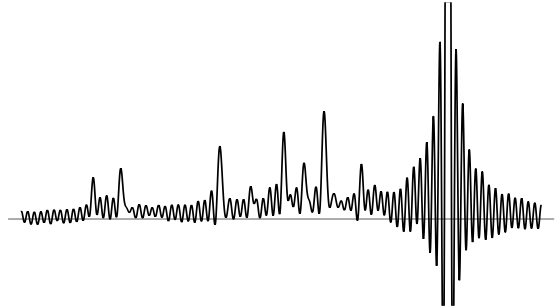,width=0.5\textwidth}
\caption{\label{zoom6} Detailed view of the shoulder to the left of
the peak at $k=1$, calculated from the sixth iteration of the
substitution rule, with the large peak intensity truncated.}
\end{figure}

Observation~\ref{mainconj} is also supported by numerical computations
of $\widehat{\gamma}$. A naive numerical analysis of a large finite
portion of $\gL$ would yield no relevant result due to the nature of
the pinwheel pattern. In particular, the number of orientations grows
only logarithmically with the radius, while the number of tiles grows
quadratically. The results above allow a more meaningful
computation. In particular, by employing the $\Sb^{1}$-symmetry, the
problem becomes one-dimensional, and by \eqref{eq:sumbess}, the
Fourier transform is expressed explicitly as a sum of weighted Bessel
functions.  As before, this approach disregards positivity of the
intensity function, and the result displayed in Figure~\ref{pin5}
shows strong oscillations with significant overshooting, similar in
kind to the ones observed in Figure~\ref{sqnum}. No smoothing of any
kind was used.

A comparison of the two diffraction images (Figures \ref{sqnum} versus
\ref{pin5}) also supports our claim about the possible radii of
pinwheel diffraction rings. One noticeable difference is the higher
peak at $k=\sqrt{5}$, which is due to the pairs of points in $\gL$
with distance $1/\sqrt{5}$. The existence of positive shoulders in
Figure~\ref{pin5} (such as that to the left of the first peak at
$k=1$, a blow-up of which is shown in Figure~\ref{zoom6}) is another
significant difference to Figure~\ref{sqnum}, and is one of the
reasons why we expect a non-vanishing radially continuous
contribution, hence giving an absolutely continuous component to
$\widehat{\gamma}$.

\section*{Acknowledgements}
It is a pleasure to thank Friedrich G\"otze, Robert V.\ Moody, Thomas
Proffen and Anthony Quas for helpful discussions. This work was
supported by the German Research Council (DFG) within the
Collaborative Research Center 701, and by EPSRC via Grant EP/D058465.

\bigskip

\parindent 0pt

\end{document}